\documentclass[brochure,12pt]{bourbaki}
\usepackage[matrix,arrow]{xy}
\usepackage{amssymb,amsfonts,amsmath,footnote}
\usepackage[francais]{babel}
\usepackage[latin1]{inputenc}
\usepackage[dvips]{graphicx}
\date{Juin 2007}
\bbkannee{59ème année, 2006-2007}
\bbknumero{980}
\title{La conjecture de Kashiwara--Vergne}
\subtitle{d'après Alekseev et Meinrenken  }
\author{Charles TOROSSIAN}
\address{Université Denis Diderot -- Paris 7\\
CNRS\\Institut de Mathématiques de Jussieu\\
175, rue du Chevaleret\\
F-75013 Paris}
\email{torossian@math.jussieu.fr}

\def\d{\mathrm{d}}

\def\g{\mathfrak g}
\def\G{\Gamma}

\def\ov{\overline}
\def\tr{\mathrm{tr}}
\def\ad{\mathrm{ad}}

\newcommand{\R}{\mathbf{R}}

\newcommand{\noi}{\noindent}

\begin{document}
\maketitle

\noindent{\bf INTRODUCTION}

En 1978, M. Kashiwara et M. Vergne   ont conjecturé dans  \cite{KV} une propriété remarquable et universelle sur  la série de Campbell-Hausdorff d'une algèbre de Lie réelle $\g$ de dimension finie. Cette propriété conjecturale admet comme corollaire l'isomorphisme de Duflo entre le centre de l'algèbre enveloppante de~$\g$ et les invariants de l'algèbre symétrique.
Cette conjecture a été démontrée en toute
généralité par A.~Alekseev et E.~Meinrenken en 2005 et publiée en 2006 à Inventiones \cite{AM}.\\

Ce texte se décompose de la façon suivante;  on rappelle dans un premier temps des résultats élémentaires sur la formule de Campbell-Hausdorff et la symétrisation. On introduit ensuite la conjecture de Kashiwara--Vergne, en expliquant ses origines et ses conséquences. La troisième section est consacrée à la preuve d'Alekseev et Meinrenken. On a résumé, dans l'appendice, la construction de Kontsevich pour la quantification des crochets de Lie qu'il nous a semblé nécessaire de rappeler pour une bonne compréhension du texte.\\

Je remercie A. Alekseev, B. Keller, D. Manchon, F. Rouvière  et M. Vergne pour leurs commentaires, suggestions et améliorations lors de la relecture de ce texte.
\section{Formule de Campbell-Hausdorff et symétrisation}

\subsection{La formule de Campbell-Hausdorff}
Soit $\g$ une algèbre de Lie de dimension finie sur  $\mathbf{R}$.
D'après le  théorème de Lie, il existe un groupe de Lie réel  $G$,
connexe et simplement connexe d'algèbre de Lie $\g$ et  une application exponentielle
notée $\exp_{\g}$ qui définit un difféomorphisme local en $0\in \g$
sur $G$.

Il en résulte que l'on peut lire la loi de groupe de $G$
en coordonnées exponentielles. C'est la fameuse formule de
Campbell-Hausdorff. En d'autre termes, pour $X, Y$
proches de $0$ dans $\g$, il existe une série en des
polynômes de Lie,  convergente et à valeurs dans $\g$, notée $Z(X, Y)$ telle que
l'on ait
$$\exp_{\g}(X)\cdot_G\exp_{\g}(Y)=\exp_{\g}\big( Z(X,Y)\big).$$
Les premiers termes de la série de Campbell-Hausdorff sont bien connus et s'écrivent

\begin{multline}
Z(X,Y)=X+Y +\frac{1}{2} [X,Y] + \frac{1}{12}[X,[X,Y]]
+\frac{1}{12}[Y,[Y,X]]+\\
\frac{1}{48}[Y,[X,[Y,X]]]-\frac{1}{48}[X,[Y,[X,Y]]]+ \cdots.
\end{multline}

Il existe de nombreuses expressions de $Z(X, Y)$ en terme de crochets itérés (cf. \S~\ref{dynkin}) ou écrites de manière récursive (cf. \cite{Va} page 118). Une difficulté  majeure concernant  la
formule de Campbell-Hausdorff est qu'il n'existe pas de base de l'algèbre de Lie
libre qui soit particulièrement commode pour  effectuer des
calculs\footnote{Les bases de Hall, par exemple sont définies de
manière récursive. Par ailleurs il existe une base qui permet
d'écrire la  formule de Campbell-Hausdorff  en utilisant des combinaisons à coefficients complexes~\cite{kly}.}. La quantification de Kontsevich donne une autre fa\c
con d'écrire la formule de Campbell-Hausdorff comme rappelé en \S~\ref{sectionKVnew}.

\subsection{Symétrisation et application exponentielle }

La symétrisation $\beta$ est un
isomorphisme d'espaces vectoriels entre l'algèbre
symétrique de $\g$ notée $S[\g]$ et l'algèbre enveloppante de $\g$ notée $U(\g)$ ; c'est une version du théorème de
Poincaré-Birkhoff-Witt. La symétrisation  commute à l'action adjointe (resp. aux dérivations $\ad X$) et
vérifie la condition pour $X\in \g$,
$$\beta(X^n)=X^n.$$On en déduit la formule, pour $X_i\in \g$,
\begin{equation}\nonumber
\beta(X_1\cdots X_n)=\frac{1}{n!}\sum\limits_{\sigma\in \Sigma_n}
X_{\sigma(1)}\cdots X_{\sigma(n)}.
\end{equation}

Il existe plusieurs fa\c cons de voir l'algèbre $S[\g]$ ; comme
algèbre symétrique, comme algèbre des fonctions polynomiales sur $\g^*$,
comme algèbre des opérateurs différentiels invariants par translation  sur $\g$ (ie. à coefficients constants)
et enfin comme algèbre pour la convolution des distributions de
support $0$. On peut voir $U(\g)$ comme l'algèbre
enveloppante universelle, l'algèbre des opérateurs
différentiels invariants à gauche sur $G$ et enfin l'algèbre
pour la convolution des distributions supportées par l'origine de $G$. \\

La formule de Taylor énonce que l'on a l'égalité de distributions
formelles $e^X=\delta_X$, avec $\delta_X$ la masse de Dirac au point
$X$. On a donc dans une complétion adéquate de~$U(\g)$ :
\begin{equation}
\beta(\delta_X)=\beta(e^X)=\sum\limits_{n\geq 0}
\frac{X^n}{n!}=\delta_{\exp_{\g}(X)},\end{equation}
o\`u $\delta_{\exp_{\g}(X)}$ est la distribution ponctuelle au point
$\exp_{\g}(X)$ dans  $G$.\par\bigskip

\noi La symétrisation envoie donc la distribution $\delta_X$ sur la
distribution $\delta_{\exp_{\g}(X)}$ ;  c'est donc l'application
exponentielle au niveau des distributions, car on a $(\exp_\g)_*
(\delta_X)=\delta_{\exp_\g(X)}$.\par\bigskip

On peut ramener, via l'application $\beta$,  le produit de $U(\g)$,
en un produit associatif dans $S[\g]$. C'est l'étoile produit de
Gutt\footnote{C'est à dire que l'étoile produit s'exprime comme un série formelle d'opérateurs bi-différentiels sur $\g^*$ à coefficients polynomiaux.} \cite{gutt}, réalisant la quantification par déformation de l'algèbre de Poisson $S[\g]$. On a donc pour $w,v$ dans $S[\g]$
\begin{equation}
w\underset{Gutt}\star v :
=\beta^{-1}(\beta(w)\beta(v)).
\end{equation}
Comme distribution de support $0$ on aura donc

\begin{lemm}\label{gutt}
Pour  $w,v$ des éléments de $S[\g]$, $\beta^{-1}(\beta(w)\beta(v))$
correspond à la distribution de support $0\in \g$ définie  pour $f$,
fonction test sur $\g$, par la formule :
$$\langle \, \beta^{-1}(\beta(w)\beta(v)) , f\, \rangle:= \langle \, w(X)\otimes v(Y), f(Z(X,Y))\, \rangle.$$
\end{lemm}

\subsection{Le centre de l'algèbre enveloppante et l'isomorphisme de
Duflo}\label{sectionisoduflo}

Comme  la symétrisation $\beta$ commute aux dérivations, c'est aussi  un
isomorphisme d'espaces vectoriels  de $S[\g]^\g$ sur $U(\g)^\g$ (les invariants pour l'action adjointe).\\

Dans  \cite{Duflo77} en utilisant les  idéaux primitifs dans l'algèbre enveloppante,
 Duflo montre que pour \textit{toute algèbre de Lie} de dimension finie sur un corps de caractéristique nulle,
 $S[\g]^\g$ et $U(\g)^\g$ sont isomorphes \textit{comme algèbres} et
exhibe un isomorphisme. Ce résultat généralise celui de Dixmier \cite{Dix}
 dans le cas nilpotent, de Duflo dans le cas résoluble \cite{Duflo70} et d'Harish-Chandra \cite{HC} dans le cas semi-simple et s'inscrit dans l'esprit de la méthode des  orbites initiée par Kirillov.

\par\bigskip

On notera par $j(X)$ le déterminant jacobien de la
fonction $\exp_{\g}$ (cf.   \S~\ref{sectiondiffexp}), à savoir la fonction définie par
\begin{equation}\label{remplacerJ}
j(X)=\mathrm{det}_{\g}\left(\frac{1-e^{-\ad X}}{\ad
X}\right)=\exp\left(-\tr_{\g}\frac{\ad X}{2}\right)
\mathrm{det}_{\g}\left(\frac{\sinh(\ad \frac{X}{2})}{\frac{\ad
X}{2}}\right).
\end{equation}
Cette fonction va intervenir de manière cruciale dans la suite.\\

Notons $S[[\g^\star]]$ l'algèbre des séries formelles en les
éléments de $\g^*$. La série formelle~$j^{\frac 12}$ est donc dans
$S[[\g^\star]]$. C'est donc un opérateur différentiel sur $\g^*$ d'ordre
infini à coefficients constants  que l'on note $j^{\frac
12}(\partial)$. La formule de Duflo s'écrit alors pour $P\in S[\g]$,
\begin{eqnarray}
\gamma(P) : =\beta\left(j^{\frac 12}(\partial) P\right).
\end{eqnarray}
C'est clairement un isomorphisme d'espaces vectoriels de $S[\g]$ sur
$U(\g)$, qui commute aux dérivations $\ad X$.

\begin{theo}[\cite{Duflo77}]\label{theoduflo}
L'application $\gamma$ ci-dessus est un isomorphisme d'algèbres de
 $S[\g]^\g$ sur $U(\g)^\g$.
 \end{theo}
Ce
théorème est  non trivial. Dans  \cite{Kont} Kontsevich montre par un argument d'homotopie que $S[\g]^\g$ et $U(\g)^\g$ sont isomorphes \textit{comme algèbres} et en déduit qu'il s'agit de l'isomorphisme de Duflo\footnote{L'argument utilise la forme a priori de l'isomorphisme obtenu comparé à celui de Duflo.}; ce résultat s'étend automatiquement aux super-algèbres de Lie et  Kontsevich montre que les algèbres de cohomologie $H(\g, S[\g])$ et $H(\g, U(\g))$  sont isomorphes\footnote{En degré $0$ on retrouve les algèbres d'invariants $S[\g]^\g$ et $U(\g)^\g$.}. Dans \cite{PT} on vérifie qu'il s'agit encore de la formule de Duflo étendue à la cohomologie\footnote{Ce résultat est aussi cité  dans \cite{shoi} comme un travail en commun avec Kontsevich, mais  non publié.}. On dispose donc d'une démonstration qui n'utilise pas la théorie des représentations.\\

Afin d'étendre l'isomorphisme de Duflo aux germes de distributions invariantes,
Kashiwara et Vergne suggèrent dans  \cite{KV} une méthode basée sur une déformation
 de la formule de Campbell-Hausdorff. Ces techniques sont  connues aujourd'hui sous le vocable ``méthode de Kashiwara--Vergne". En quelque sorte on cherche à lire l'isomorphisme de Duflo sur la formule de Campbell-Hausdorff.

\section{La conjecture combinatoire de Kashiwara--Vergne}
\subsection{Notations}
Soient $(e_i)_{i=1,\ldots, d}$ une
base de $\g$, $(e_i^*)_{i=1,\ldots, d}$ la base duale et $X=\sum_{i=1}^{d} x_i e_i$. Pour  $X\mapsto A(X)$ une
fonction régulière de $\g$ dans $\g$ (c'est à dire un champ de vecteurs sur $\g$) on
désigne le champ adjoint associé par  \begin{equation}\label{notations}[X,
A(X)]\cdot\partial_{X}=\sum_{i=1}^{d} \langle \, e_i^*, [X, A(X)]\,
\rangle \frac{\partial}{\partial x_i}. \end{equation}

\noi  On note aussi $\partial_X A$ la différentielle
de $A$ en $X$ ; c'est une application linéaire de $\g$ dans $\g$. \\

\subsection{Énoncé de la conjecture combinatoire}

Soit $\g$ une algèbre de Lie de dimension finie sur $\R$. Pour $X, Y \in\g$ on note  $Z(X,Y)$ la série de Campbell-Hausdorff
définie par \[Z(X,Y)=\log\left(\exp_\g(X)\underset{G}\cdot\exp_\g(Y)\right).\]
 Dans tout ce qui suit nous travaillons
au niveau des séries formelles mais des arguments élémentaires
montrent que toutes les séries formelles que nous manipulons sont
convergentes dans un voisinage de $(0,0)$.

La conjecture combinatoire de Kashiwara--Vergne \cite{KV} s'énonce de
la manière suivante :

\begin{theo}[Conjecture KV 78]\label{conjKV}Notons $Z(X, Y)$ la série de Campbell-Hausdorff. Il existe des séries $F(X,Y)$ et $G(X,Y)$ sur $\g\oplus\g$ sans
terme constant et à valeurs dans $\g$ telles que l'on ait
\begin{equation}\label{KV1}
X+Y-\log(\exp_\g(Y)\cdot \exp_\g (X)) =\big(1-e^{-\ad X}\big)F(X,Y)+ \big(e^{\ad Y}
-1\big)G(X,Y)
\end{equation}
et telle que l'identité de trace suivante soit vérifiée
\begin{equation}\label{KV2}\tr_{\g}(\ad X\circ \partial_XF+\ad Y\circ \partial_YG)=
\frac{1}{2}\tr_{\g} \left(\frac{\ad X}{e^{\ad
X}-1}+\frac{\ad Y}{e^{\ad Y}-1} -\frac{\ad Z(X, Y)}{e^{\ad
Z(X, Y)}-1}-1\right).
\end{equation}
\end{theo}
\begin{rema} La conjecture porte sur l'existence d'un couple $(F, G)$ de solutions universelles\footnote{C'est à dire des éléments d'une complétion de l'algèbre de Lie libre universelle engendrée par $X, Y$. L'algèbre de Lie libre
sur un espace vectoriel $V$ est
naturellement graduée ainsi que son algèbre enveloppante,
qui est l'algèbre associative libre sur $V$;  par exemple, l'élément $[X,Y]=XY-YX$ est de
degré $2$. La complétion consiste à considérer
les séries de termes homogènes dont les degrés tendent vers l'infini \cite{Serre}.} vérifiant les équations ci-dessus.  Par ailleurs, si un tel couple convient  alors le couple $$\left(G(-Y, -X), F(-Y, -X)\right)$$ est aussi une solution. On peut donc rechercher des solutions symétriques c'est à dire vérifiant $G(X, Y)=F(-Y, -X)$. Pour de telles solutions on s'aperçoit facilement que les termes à l'ordre $1$  en $Y$ sont uniquement déterminés \cite{AP}, ce qui peut faire espérer l'unicité d'une solution symétrique\footnote{Ce point est encore conjectural.}.\end{rema}

\begin{rema} L'équation (\ref{KV1}) peut se résoudre complètement dans l'algèbre tensorielle engendrée par $X, Y$.  Puis  en utilisant l'idempotent de Dynkin (cf. \S~\ref{dynkin}) on peut exhiber   toutes les solutions de (\ref{KV1}) (cf. \cite{bur}). La difficulté majeure de cette conjecture est donc l'équation de trace (\ref{KV2}).
\end{rema}

\begin{rema}  Les équations (\ref{KV1}) et (\ref{KV2}) forment un système affine à coefficients rationnels. Si on
dispose d'une solution à coefficients réels, alors il existera une solution à coefficients rationnels. Expliciter une solution rationnelle est un problème intéressant que l'on peut poser.
\end{rema}
\subsubsection{La solution conjecturale de Kashiwara--Vergne} Dans leur article, Kashiwara et Vergne proposent un couple symétrique  $(F^0, G^0)$ de séries de Lie universelles et montrent, \textit{dans le cas résoluble},  que ce couple vérifie la condition de trace (\ref{KV2}). Nous suivons l'article \cite{Rou86} qui propose une réécriture de ce couple conjectural.\\

Notons  $\psi$ la  fonction analytique au voisinage de $0$  définie par  $$\psi(z)= \frac{e^z-1-z}{(e^z-1)(1-e^{-z})}.$$
Soit $Z(t)=Z(tX, tY)$ et posons  $$F^1(X, Y)= \left(\int_0^1 \frac{1-e^{-t \,\ad X}}{1-e^{-\ad X}} \circ \psi(\ad Z(t)) \d t\right) (X+Y)$$
et $G^1(X, Y)=F^1(-Y, -X)$. Posons alors $$F^0(X, Y)=\frac12 \left(F^1(X, Y) +e^{\ad X} F^1(-X, -Y)\right) +\frac14 \left(Z(X, Y) -X\right)$$ et  $G^0(X, Y)=F^0(-Y, -X)$.\footnote{Écrivons  pour  simplifier $x=\ad X$ et $y=\ad Y$. Des calculs  fastidieux dans les années $80$, mais que l'on peut maintenant effectuer sur ordinateur montrent que les  premiers termes s'écrivent (jusqu'à l'ordre $4$) : $F^0(X, Y)= \frac14 Y + \frac1{24}xY -\frac1{48}x^2Y-\frac1{48}yxY-\frac1{180}x^3Y-\frac1{480}yx^2Y+\frac1{360}y^2xY + \ldots$}\\

Par construction (cf. \cite{KV}, \cite{Rou86}) ce couple $(F^0, G^0)$ vérifie la première équation~(\ref{KV1}). Dans leur article Kashiwara et Vergne  conjecturent que  $(F^0, G^0)$ vérifie l'équation de trace~(\ref{KV2}) et  vérifie ce fait  dans le cas des algèbres de Lie résolubles.\\

\subsubsection{Historique des résultats} Dans \cite{rou81}, Rouvière  vérifient la conjecture  dans le cas $\mathrm{sl}(2, \R)$  pour le couple $(F^0, G^0)$ ci-dessus.\\

En 1999  dans \cite{Ve} Vergne démontre la conjecture
 dans le cas quadratique, c'est à dire pour une algèbre de Lie munie d'une forme bilinéaire
invariante et non dégénérée;  les algèbres réductives mais aussi $\g\oplus \g^*$ avec $\g$ quelconque, sont quadratiques. L'article \cite{Ve}  suit les idées de \cite{AM0} concernant l'isomorphisme de Duflo pour les algèbres de Lie quadratiques. Dans \cite{AM2} Alekseev et Meinrenken proposent,  toujours dans le cas quadratique, une solution différente en utilisant la géométrie de Poisson et le ``Moser trick". Toutefois dans  \cite{AP} il est montré que ces solutions du cas quadratique ne sont pas universelles, c'est à dire qu'elles ne résolvent  pas la conjecture pour toutes les algèbres de Lie.\\

Enfin  A. Alekseev et E. Meinrenken \cite{AM} ont démontré
en mai 2005 la conjecture combinatoire  en utilisant une déformation  de la série de Campbell-Hausdorff décrite dans \cite{To} et résultant
de la quantification de Kontsevich \cite{Kont}. Le lien avec le couple conjectural $(F^0, G^0)$ n'y est cependant pas abordé.

\subsection{Origine et conséquences de cette conjecture} Cette égalité sur les traces peut sembler étrange mais
elle est une consé\-quence naturelle de l'intégration par partie.
Expliquons un peu tout ceci ce qui motivera le lecteur.

\subsubsection{Transport de la convolution}

Un des buts de l'article de Kashiwara--Vergne est de démontrer, pour les distributions invariantes, le transport du produit de convolution par l'application exponentielle.\\

Plus précisément il s'agissait de montrer le résultat conjectural \footnote{Conjecture aussi formulée par Raïs.} suivant assurant que l'on peut transporter la convolution sur $\g$ en la convolution sur $G$. Ce point est important car  les caractères des représentations irréductibles sont des solutions propres pour le centre de $U(\g)$ et exprimés en coordonnées exponentielles deviennent alors des solutions propres des opérateurs différentiels invariants à coefficients constants\footnote{L'évaluation sur une orbite définit alors un caractère pour $S[\g]^\g$.}.\\

Ce théorème  fut démontré par la suite dans \cite{ADS}, \cite{AST}, \cite{Mo} comme corollaire de la quantification de Kontsevich.

\begin{theo}[\cite{ADS}, \cite{AST}, \cite{Mo}]\label{KVR} Soient $w$ et $v$
deux germes de distributions invariantes au voisinage de $0$ dans $\g$ et
vérifiant une certaine condition de support \footnote{On peut demander par exemple que les supports asymptotiques de $u$ en $0$  (resp. $v$),  noté $C_u$ (resp. $C_v$), vérifient $C_u\cap -C_v=\{0\}$. } afin d'assurer un sens à
la convolution. On a
\begin{equation}\label{eqKVbis}
\langle \, \, w\otimes v\,, \,\frac{j^{\frac 12}(X)j^{\frac
12}(Y)}{j^{\frac 12}(Z(X,Y))}f(Z(X,Y))\,\, \rangle=\langle \,\,
w\otimes v\,, \,f(X+Y)\,\, \rangle
\end{equation}
avec  $f$ une   fonction $C^{\infty}$ dans un voisinage de $0$ et à support compact.

\end{theo}

\begin{defi} On appellera fonction de densité le quotient  $$D(X, Y) := \frac{j^{\frac 12}(X)j^{\frac
12}(Y)}{j^{\frac 12}(Z(X,Y))}.$$ \end{defi}
\subsubsection{Déformation par dilatation}
L'idée de base de l'article \cite{KV} est de considérer la
déformation naturelle de  l'algèbre de Lie $\g$ qui consiste à
remplacer le crochet $[X,Y]$ par  $t[X,Y]$ pour $t\in[0,1]$. La série de
Campbell-Hausdorff est changée en
$Z_{t}(X,Y)=\frac{1}{t}Z(tX,tY)$. Remarquons que cette expression est bien définie pour $t=0$ et on obtient
$Z_0(X, Y)=X+Y$.\\

On déduit de la différentielle de l'application
exponentielle (cf. \S \ref{equivalence}), que l'équation~(\ref{KV1}) est équivalente à l'équation différentielle suivante
\begin{equation}\label{KV1bis}
\frac{\partial}{\partial t}Z_{t}(X,Y)=[X,
F_{t}(X,Y)]\cdot\partial_{X}Z_{t}(X,Y) +[Y,
G_{t}(X,Y)]\cdot\partial_{Y}Z_{t}(X,Y),
\end{equation}
o\`u on a noté $F_{t}(X,Y)=\frac{1}{t}F(tX,tY)$\footnote{Comme $F(0, 0)=(0, 0)$ cette expression est bien définie pour $t=0$.}.\\

Notons  $q(X)$ la fonction
$\det_{\g}\left(\frac{\sinh \frac{\ad X}2}{\frac{\ad X}2}\right)$. En utilisant  la formule\footnote{Voir \S\ref{sectiondiffexp} pour les nombres de Bernoulli $b_n$.}
$$\ln(q(X))=\sum_{n\geq 1} \frac{b_{2n} \tr_\g (\ad X)^{2n}}{(2n)! 2n},$$ on trouve facilement \cite{KV} :
\begin{equation}\label{eqdiffjbis}
j^{-\frac 12}(tX)\frac{\partial}{\partial t}j^{\frac
12}(tX)=\frac{1}{2}\tr_{\g}\left(\frac{\ad X} {\exp(t\,
\ad X)-1}-\frac{1}{t}\right).
\end{equation}

\noi

\noi
Compte tenu de la déformation en le paramètre $t$, il suffit de démontrer  qu'on a  l'égalité pour tout $t\in[0,1]$,

\begin{equation}\label{eqKVtbis}
\langle \,\, w\otimes v\,, \,\frac{j^{\frac 12}(t\,X)j^{\frac
12}(t\,Y)}{j^{\frac 12}(t\,Z_{t}(X,Y))}f(Z_{t}(X,Y))\,\, \rangle=
\langle \,\, w\otimes v\,, \,f(X+Y)\,\, \rangle.
\end{equation}
L'idée est maintenant simple, il suffit de demander  que la dépendance en
$t$ soit triviale, c'est à dire que  la dérivée par rapport à $t$ soit
nulle.\\

\subsubsection{Calcul de la dérivée en $t$} Notons $D_t(X,Y)$
la fonction de densité
$$D_t(X,Y)=D(tX, tY)=\frac{j^{\frac 12}(tX)j^{\frac 12}(tY)}{j^{\frac 12}(tZ_{t}(X,Y))}.$$
 On peut  remplacer $j$ par $q$ dans cette formule, la fonction de densité reste la même.
On a facilement compte tenu des équations (\ref{KV1bis}) et
(\ref{eqdiffjbis}),

\begin{multline}\label{deriveeA}\frac{\partial}{\partial t}D_t(X,Y)=
\frac{1}{2}\tr_{\g}\left(\frac{\ad X}
{e^{t\,\ad X}-1}+\frac{\ad Y}{e^{t\,\ad Y}-1}
-\frac{\mathrm{ad}Z_t(X,Y)}{e^{t\,\mathrm{ad}Z_{t}(X,Y)}-1}-
\frac{1}{t}\right)D_t(X,Y)+\\
[X, F_{t}(X,Y))]\cdot\partial_{X}D_t(X,Y) +[Y,
G_{t}(X,Y))]\cdot\partial_{Y}D_t(X,Y).
\end{multline}
Pour simplifier on va noter :
\begin{equation}\label{T(X,Y)}
T(X,Y)=\frac{1}{2}\tr_{\g}\left(\frac{\ad X}
{e^{\ad X}-1}+\frac{\ad Y}{e^{\ad Y}-1}
-\frac{\mathrm{ad}Z(X,Y)}{e^{\mathrm{ad}Z(X,Y)}-1}-1\right).
\end{equation}
Par conséquent le premier  terme du second membre de (\ref{deriveeA}) est
$\frac{1}{t}T(tX,tY)$. Ce calcul se justifie comme suit; le premier terme dans le membre de droite (\ref{deriveeA}) résulte de la dérivée par rapport à $t$ dans
les termes $j^{\frac 12}(t\cdot)$ et le second terme résulte de la
dérivée en $t$ dans  $Z_{t}$. Plus précisément compte tenu de
(\ref{KV1bis}) il vient que pour toute fonction $\phi$ on a
\begin{equation}
\frac{\partial}{\partial t} \phi\left(Z_{t}(X,Y)\right)=[X,
F_{t}(X,Y)]\cdot\partial_{X}\phi\left(Z_{t}(X,Y)\right) + [Y,
G_{t}(X,Y)]\cdot\partial_{Y}\phi\left(Z_{t}(X,Y)\right).
\end{equation}

Le  champ de vecteurs $[X, F_{t}(X,Y)]\cdot\partial_{X} +[Y,
F_{t}(X,Y)]\cdot\partial_{Y}$ agit trivialement sur la fonction
$j^{\frac 12}(tX)j^{\frac 12}(tY)$ car cette dernière est invariante
en chaque variable sous l'action adjointe,n par conséquent la dérivée
du terme
en  $Z_{t}$ s'écrit bien comme annoncée.\\

\noi On peut maintenant terminer le calcul de la dérivée dans
(\ref{eqKVtbis}). Il vient
\begin{multline}\label{diffdensite}
\frac{\partial}{\partial t}\left(D_t(X,Y) f(Z_{t}(X,Y))\right)=\\
\Big([X, F_{t}(X,Y)]\cdot\partial_{X} +[Y,
G_{t}(X,Y)]\cdot\partial_{Y}\Big)\left(D_t(X,Y) f(Z_{t}(X,Y))\right) +\\
 \frac{1}{t}T(tX,tY)D_t(X,Y) f(Z_{t}(X,Y)).\end{multline}
On est donc amené  à calculer l'action à droite\footnote{En effet
les distributions sont plutôt un module à droite sur les champs
de vecteurs. Cette action à droite n'est utilisée qu'à cet endroit du texte.} du champ de vecteurs $[X,
F_{t}(X,Y)]\cdot\partial_{X} +[Y, F_{t}(X,Y)]\cdot\partial_{Y}$ sur la distribution $w\otimes v$.
Compte tenu de l'invariance de cette distribution on a,
\begin{multline}\label{invarianceUVbis}
w\otimes v\Big([X, F_{t}(X,Y))]\cdot\partial_{X} +[Y,
F_{t}(X,Y))]\cdot\partial_{Y}\Big)=\\
-w\otimes v\Big(\tr_{\g}\big(\ad X\circ
\partial_XF_{t}(X,Y)+\ad Y\circ \partial_Y G_{t}(X,Y)\big)\Big).
\end{multline}
Pour conclure au transport de la convolution dans  (\ref{eqKVbis})
 il suffit de demander que l'on ait pour tout $t$ :
\begin{equation}\label{tracebis}
\frac{1}{t}T(tX,tY)-\tr_{\g}\big(\ad X\circ
\partial_XF_{t}(X,Y)+\ad Y\circ \partial_YG_{t}(X,Y)\big)=0.
\end{equation}

\noi \textbf{La conjecture combinatoire de Kashiwara--Vergne est
précisément cette égalité.}\\
\begin{rema}\label{equivalencebis} Si l'égalité (\ref{KV1}) (ou bien (\ref{KV1bis})) est vérifiée, alors l'équation (\ref{KV2}) est équivalente
à l'équation (\ref{KV2bis}) suivante

\begin{multline}\label{KV2bis}
\frac{\partial}{\partial t}D_t(X,Y)= \Big([X,F_{t}(X,Y)]\cdot\partial_{x} +[Y,
G_{t}(X,Y)]\cdot\partial_{y}\Big)D_t(X,Y) +\\
\Big(\tr_{\g}\big(\ad X\circ \partial_XF_{t}(X,Y)+
\ad Y\circ \partial_YG_{t}(X,Y)\big)\Big)D_t(X,Y) .\end{multline}

\end{rema}
\subsubsection{Isomorphisme de Duflo} On déduit comme cas particulier du transport de la convolution (\ref{eqKVbis}) le corollaire suivant :
\begin{coro}
Le théorème \ref{KVR} généralise l'isomorphisme de Duflo.
\end{coro}
\noindent{\sc Preuve}  --- En effet pour $P,
Q\in S[\g]^\g$ considérées comme des distributions invariantes de
support $0$ et  pour toute fonction test $f$, on a d'après le théorème \ref{KVR} appliqué à $j^{1/2} f$,
\begin{equation}\nonumber \langle \, P\otimes Q,
j^{\frac 12}(X)j^{\frac 12}(Y)f(Z(X,Y))\, \rangle=\langle \, P\otimes Q, (j^{1/2}f)(X+Y)\,
\rangle=\langle \, j^{\frac 12}(\partial) (PQ), f \, \rangle.
\end{equation}
On a donc
\begin{equation}\nonumber \langle \, j^{\frac 12}(\partial) P\otimes  j^{\frac 12}(\partial) Q,
f(Z(X,Y))\, \rangle=\langle \, j^{\frac 12}(\partial)(PQ), f \,
\rangle.
\end{equation}
Et d'après le lemme   \ref{gutt},  le membre de gauche correspond à l'élément  \[\beta^{-1}\left(\beta\left(j^{\frac
12}(\partial) P\right)\beta\left(j^{\frac 12}(\partial) Q\right)\right), \] tandis que
le membre de droite correspond à  $j^{\frac 12}(\partial) (PQ)$. On en déduit que l'on a   bien la formule de Duflo sur les éléments invariants,

\[\beta\Big(j^{\frac
12}(\partial) (PQ)\Big)=\beta\left(j^{\frac 12}(\partial)
P\right)\beta\left(j^{\frac 12}(\partial) Q\right).\]

\section{Preuve de la conjecture combinatoire de Kashiwara--Vergne}
On explique dans cette section les deux résultats principaux de \cite{AM}, à savoir la preuve de la conjecture de Kashiwara--Vergne et l'extension de l'isomorphisme de Duflo à toute la cohomologie comme corollaire. Le texte suit essentiellement l'article \cite{AM}.\\

On pose dans un premier temps des définitions utiles.
\subsection{Notations}

\subsubsection{Application de Duflo sur les distributions}
On fixe un voisinage symétrique $\mathcal{U}$ de $0$ dans $\g$ sur lequel $\exp_\g$ est un difféomorphisme sur son image $U=\exp_\g(\mathcal{U})$. On suppose  $j>0$ sur $\mathcal{U}$.  L'application de Duflo sur les distributions à support compact
$$\mathrm{Duf}=(\exp_\g)_\star \circ j^{1/2}: \mathcal{D}'_{comp}(\mathcal{U}) \quad \rightarrow \quad \mathcal{D}_{comp}'(U),$$ est alors bien définie. On rappelle que l'on a noté $D(X, Y)=\frac{j^{1/2}(X)j^{1/2}(Y)}{j^{1/2}(Z(X, Y))}$ la fonction de densité.\\

\subsubsection{Produit $m$}
On fixe  alors un voisinage $\mathcal{O}$ de $0$ dans $\g$ suffisamment petit, tel que $\exp_\g(\mathcal{O})\exp_\g(\mathcal{O})\subset \mathcal{U}$ et tel que la série de Campbell-Hausdorff $Z$ définisse une fonction  régulière de $\mathcal{O}\times \mathcal{O}$ dans $\mathcal{U}$.\\

La convolution des distributions à support compact dans $\exp_\g(\mathcal{O})$ est alors bien définie et on peut remonter le  produit de convolution de $G$ via l'application de Duflo en un produit sur $ \mathcal{D}'_{comp}(\mathcal{O})$ ; c'est la formule du théorème~\ref{KVR}. On note ce produit

\begin{equation}\label{defim}m=Z_\star \circ D : \quad  \mathcal{D}'_{comp}(\mathcal{O}\times \mathcal{O})\quad  \rightarrow \quad \mathcal{D}'_{comp}(\mathcal{U}).\end{equation}

\subsubsection{Déformation par dilatation $m_t$}

Utilisons la déformation du crochet de Lie : $[X, Y]_t= t [X, Y]$ et notons $\g_t$ cette nouvelle algèbre de Lie. Pour $t=0$ on trouve le crochet trivial, et pour $t\neq 0$ c'est une algèbre isomorphe à $\g$.\\

La série de Campbell-Hausdorff pour $\g_t$ s'écrit $Z_t(X, Y)=\frac1tZ(tX, tY)$ et la fonction de densité vaut $D_t(X, Y)=D(tX, tY)$. On en déduit comme précédemment un produit $m_t$ sur $\mathcal{D}'_{comp}(\frac1t\mathcal{O})$

 $$m_t=(Z_t)_\star \circ D_t : \quad  \mathcal{D}'_{comp}(\frac1t\mathcal{O}\times \frac1t\mathcal{O})\quad  \rightarrow \quad \mathcal{D}'_{comp}(\frac1t\mathcal{U}).$$

\subsubsection{Dérivée de Lie} La dérivée de Lie de l'action adjointe sur les fonctions $f\in\mathcal{C}^\infty(\g) $ est définie pour $v\in \g$ par $$(L(v) f )(X)= [X, v]\cdot \partial_X f.$$

Si $u$ est une distribution et $f$ une fonction test, alors par définition\footnote{Il n'y aurait pas de signe $-$ si on prenait l'action à droite.}  on a $\langle L(e_i) u, f\rangle= -\langle u, L(e_i) f \rangle$. Une distribution est dite invariante sur elle est annulée par les $L(e_i)$ pour $i=1, \ldots, n$\footnote{On rappelle que $(e_i)_i$ est une base de $\g$.}.\\

Soit $X\mapsto A(X)=A_1(X) e_1 +\ldots + A_n(X) e_n$ est un champ de vecteurs sur $\g$. Alors la dérivée de Lie $L(A)$ vaut $\sum_i A_i L(e_i)$ sur les fonctions et $\sum_i L(e_i) \circ A_i$ sur les distributions.\\

\subsubsection{Structure de Lie sur \; $\mathcal{C}^\infty(\g\times \g, \g\times \g)$}\label{crochetLie}
On définit une structure de Lie $[\quad , \quad ]_{Lie}$  sur les fonctions $\mathcal{C}^\infty(\g\times \g, \g\times \g)$ en demandant que
$$\beta= (F(X, Y), G(X, Y)) \quad \mapsto \quad \Xi^\beta= [X, F(X, Y)]\cdot \partial_X + [Y, G(X, Y)]\cdot \partial _Y$$
soit un morphisme lorsqu'on munit les champs de vecteurs sur $\g\times\g$ du crochet standard.\\

Explicitement on a pour $\beta, \gamma \in \mathcal{C}^\infty(\g\times \g, \g\times \g)$ en utilisant la dérivée de Lie $L(\beta)$ du champ de vecteurs sur $\g\times\g$\footnote{La dérivée de Lie est pour l'action adjointe de $\g\times\g$ sur les fonctions.} ,

$$[\beta, \gamma]_{Lie}= L(\beta) \gamma -L(\gamma) \beta + [\beta, \gamma]_{\g\times\g}.$$

\subsection{Reformulation de la conjecture}

Dans leur article Alekseev et Meinrenken utilisent une reformulation plus  géométrique de la conjecture de Kashiwara--Vergne.

Pour cela il est commode d'utiliser la base de $\g\times \g$ et les dérivées de Lie de l'action adjointe de $\g\times\g$. On notera $(\hat{e_i})_{i=1, \ldots, 2n}$ une base de $\g\times \g$ (on prend deux copies de la base $e_i$) et  $L(\hat{e_i})$ la dérivée de Lie de l'action de $\g\times \g$ sur $\mathcal{D}'_{comp}(\mathcal{O}\times \mathcal{O}).$\\

\subsubsection{Deux reformulations équivalentes de la conjecture de Kashiwara--Vergne}
Soit $\beta \in  \mathcal{C}^\infty(\mathcal{O}\times \mathcal{O}, \g\times \g)$, tel que $\beta(0, 0)=0$. On suppose que $\beta$ est  un couple vérifiant la conjecture de Kashiwara--Vergne. Utilisons la base ci-dessus et écrivons $$\beta=(F(X, Y), G(X, Y))= \sum_{1\leq i \leq 2n} \beta^i \hat{e_i}$$
et $\beta_t=\frac1t \beta(tX, tY)= (F_t(X, Y), G_t(X, Y))= \sum\limits_{1\leq i \leq 2n} \beta^i_t \hat{e_i}$.\\

\noi \textbf{Première reformulation:} Les équations (\ref{KV1}) et (\ref{KV2}) sont équivalentes aux équations (\ref{KV1bis}) et (\ref{KV2bis}) ( cf. \S~\ref{equivalence} et la remarque \ref{equivalencebis} ) que l'on peut écrire de manière plus condensée en utilisant les dérivées de Lie sur les fonctions :

\begin{equation} \partial_t Z_t= \left(\sum_{1\leq i \leq 2n} \beta_t^i L(\hat{e_i})\right)Z_t\end{equation}et
\begin{equation} \partial_t D_t=  \sum_{1\leq i \leq 2n} L(\hat{e_i})(\beta_t^i D_t)= \left(\sum_{1\leq i \leq 2n}  L(\hat{e_i})\circ \beta_t^i\right) D_t.\end{equation}\\

\noi \textbf{Deuxième reformulation:} L'idée astucieuse de \cite{AM} est de constater que ces équations peuvent se condenser en une seule. La preuve est élémentaire et consiste à vérifier l'assertion sur les distributions de dirac $\delta_{p, q}$ au point $(p, q)\in \g\times \g$.

\begin{prop}\label{KVI} Les équations (\ref{KV1}) et (\ref{KV2}) sont équivalentes à l'équation suivante où la dérivée de Lie porte sur les distributions
\begin{equation} \partial_t m_t=- m_t \circ \left(\sum_{1\leq i \leq 2n} \beta_t^i L(\hat{e_i})\right).
\end{equation}
\end{prop}

\noi \textbf{Notation :} On notera $V(\beta)$ l'opérateur agissant sur les distributions $\sum\limits_{1\leq i \leq 2n} \beta_t^i L(\hat{e_i})$.\\

La conjecture de Kashiwara--Vergne s'écrit donc plus simplement $$\partial_t m_t=- m_t \circ V(\beta_t).$$

\begin{rema} Comme on le constate ce n'est pas la dérivée de Lie $L(\beta_t)$  sur les distributions qui intervient mais l'opérateur sur les distributions $V(\beta_t)= \sum_{1\leq i \leq 2n} \beta_t^i L(\hat{e_i})$. Toutefois on vérifie (cf. \cite{AM})  que l'on a la relation importante suivante pour $\beta, \gamma \in \mathcal{C}^\infty(\mathcal{O}\times \mathcal{O}, \g\times \g)$,

$$V([\beta, \gamma]_{Lie})=[V(\beta), V(\gamma)],$$c'est à dire $V$ est un homomorphisme de Lie.

\end{rema}
\subsection{Equation de courbure nulle}
 Dans \cite{To} on définit, grace à  la quantification de Kontsevich (cf. Appendice A \S~\ref{sectiondeformationkont}), une déformation $Z_u(X, Y)$ (resp. $D_u(X, Y)$) de la série de Campbell-Hausdorff (resp. de la fonction de densité).\\

  On écrit maintenant ces équations de déformation en termes analogues à la proposition \ref{KVI}. D'après \S~\ref{sectiondeformationkont} théorème~\ref{theoT}, il existe une fonction $\gamma_u \in \mathcal{C}^\infty(\mathcal{O}\times \mathcal{O}, \g\times \g)$ analytique en $u\in[0, 1]$, vérifiant $\gamma_u(0, 0)=(0, 0)$  et donnée par des séries de Lie universelles convergentes  et il existe  une déformation $\widehat{m_u}$ (définie comme en  (\ref{defim}) avec  $Z_u$ et $D_u$ )   du produit sur les distributions telles que l'on ait

$$\partial_u \widehat{m_u} =-\widehat{m_u} \circ V(\gamma_u).$$
 Pour $u=0$ on retrouve $m_0$  le produit standard dans $\g$ et pour $u=1$ on retrouve le produit $m$. Il est important de remarquer que la déformation $m_u$ n'est pas un produit associatif, contrairement à $m_t$. \\

 En considérant le paramètre de déformation par dilatation $t\in[0,1]$ on en déduit que l'on a aussi par dilatation

 \begin{equation}\label{eqdiffmut}\partial_u \widehat{m_{u,t}} =-\widehat{m_{u,t}} \circ V(\gamma_{u,t})\end{equation}avec $\gamma_{u,t}(X, Y)=\frac1t \gamma_u(tX, tY)$ et $\widehat{m_{u,t}}$   définie comme en  (\ref{defim})  avec  $\frac1tZ_u(tX, tY)$ et $D_u(tX, tY)$. On a $\widehat{m_{u=0,t}}=m_0$ et $\widehat{m_{u=1,t}}=m_t$.\\

 L'idée de Alekseev et Meinrenken est de construire à partir de la solution  $\gamma_u$ une solution $\beta_t$ de la proposition~\ref{KVI}. Pour cela il faut résoudre une équation de courbure.

 \begin{prop} Il existe  une série de Lie universelle $\beta_{u, t}$  à valeurs dans $\g\times \g$ convergente dans un voisinage de $(0, 0)$, telle que $\beta_{0, t}=(0, 0)$ et vérifiant l'équation
 $$\partial_u \beta_{u, t} -\partial_t\gamma_{u, t} + [\beta_{u,t}, \gamma_{u,t}]_{Lie}=0.$$ On a $\beta_{u, t}(X, Y)=\frac1t \beta_{u, t=1}(tX, tY)$.
\end{prop}

\noindent{\sc Preuve}  --- C'est une équation différentielle en $u$ linéaire  avec second membre. Il  y a donc unicité et l'assertion sur la dilatation  résulte de l'égalité  $\gamma_{u, t}(X, Y)= \frac1t \gamma_{u, t=1}(tX, tY)$. On peut résoudre formellement cette équation en utilisant les résolvantes\footnote{On pourrait utiliser aussi les séries de Magnus.}. On trouve facilement

$$\beta_{u,t}= \sum\limits_{n\geq 0} \int_{0\leq u_0 \leq u_1 \leq \ldots \leq u_n \leq u} \d u_0 \ldots \d u_n \ad \gamma_{u_n, t} \ldots \ad \gamma_{u_1, t} \partial_t \gamma_{u_0, t}.$$
A partir de cette expression on montre la convergence et d'après la définition du crochet $[\quad, \quad ]_{Lie}$ \S~\ref{crochetLie} il est clair que l'on ne manipule que des séries de Lie. \\

\begin{rema} Une autre façon de voir l'équation de courbure est la suivante. On résout formellement dans le groupe de Lie associé à la structure de Lie $[\quad, \quad ]_{Lie}$, l'équation $\partial_u g(u, t)=-g(u, t) \gamma_{u, t}$ avec condition initiale $g(0, t)=1$. On a alors facilement $$\beta_{u, t}=-g(u, t)^{-1} \partial_t g(u, t).$$
\end{rema}

\subsection{Solution d'Alekseev et Meinrenken à la conjecture de Kashiwara--Vergne}
On conclut maintenant par le théorème suivant qui résout la conjecture de Kashiwara--Vergne.
\begin{theo}[\cite{AM}] La série de Lie universelle $\beta_{u=1, t}=(F_t(X, Y), G_t(X, Y))$ résout la conjecture de Kashiwara--Vergne.
\end{theo}

\noindent{\sc Preuve}  --- Le formalisme de la remarque précédente est très pratique pour comprendre la preuve. Le point clef est que l'on  a la formule suivante, $$\widehat{m_{u,t}} = m_0 \circ V(g(u, t)).$$ En effet les deux membres coincident pour $u=0$ et vérifient la même équation différentielle (\ref{eqdiffmut}). On a en effet
\begin{multline}\partial_u \left( m_0 \circ V(g(u, t))\right) = m_0 \circ V(\partial_u g(u,t))= \\- m_0 \circ V( g(u,t) \gamma_{u,t})= -\left(m_0 \circ V(g(u, t))\right)\circ V(\gamma_{u, t})\end{multline}
car $V$ est un homomorphisme d'algèbres de Lie  pour $[\quad, \quad]_{Lie}$.\\

\noi Comme on a $\beta_{u, t}=-g(u, t)^{-1} \partial_t g(u, t)$, on en déduit $$\partial_t \widehat{m_{u,t}}= -\widehat{m_{u,t}}\circ V(\beta_{u, t}).$$Or on a  $\widehat{m_{u=1,t}}= m_t$  il vient donc l'équation $\partial_t m_t=- m_t \circ V(\beta_t)$ avec $\beta_t=\beta_{u=1,t}$. Ceci  résout la conjecture de Kashiwara--Vergne d'après la reformulation de la proposition~\ref{KVI}.

\subsection{Prolongement de l'isomorphisme de Duflo en cohomologie}
Considérons le morphisme  d'algèbres de Lie, $\psi: \g \rightarrow \g\times\g$ donné par l'injection diagonale. L'application duale $\psi^*$ s'étend aux algèbres extérieures; $$\psi^*: \bigwedge(\g^*\oplus \g^*)=\bigwedge\g^* \otimes \bigwedge \g^* \longrightarrow \bigwedge \g^*.$$ C'est le produit dans l'algèbre extérieure $\bigwedge \g^*$.\\

Si $\mathcal{M}$ est un $\g$-module, on notera $d_{\mathcal{M}}$ la différentielle de Chevalley-Eilenberg sur $\mathcal{M}\otimes \bigwedge \g^*$ et  $C(\g, \mathcal{M})$ le complexe  des cochaînes. Comme d'habitude on notera par $H(\g, \mathcal{M})$ l'algèbre de cohomologie.\\

Considérons $\mathcal{D}'_{comp}(\mathcal{O}\times\mathcal{O})$ comme un $\g\times \g$-module et $\mathcal{D}'_{comp}(\mathcal{O})$ comme un $\g$-module. Notons $M_t=m_t\otimes \psi^*$. C'est une application de complexes de
$C\left(\g\times \g, \mathcal{D}'_{comp}(\mathcal{O}\times\mathcal{O})\right)$ dans $C\left(\g, \mathcal{D}'_{comp}(\mathcal{O})\right)$. On a $$\partial_t M_t= \partial_t m_t\otimes \psi^*.$$
Dans \cite{AM} on montre que $V(\beta) \otimes \psi^*$, qui est une application de complexes de $C\left(\g\times \g, \mathcal{D}'_{comp}(\mathcal{O}\times\mathcal{O})\right)$ dans $C\left(\g, \mathcal{D}'_{comp}(\mathcal{O}\times\mathcal{O})\right)$, est homotopiquement triviale. Plus précisément on a\footnote{Ici intervient le fait que $\beta$ est une application $\g$-équivariante.}
$$V(\beta) \otimes  \psi^*=(1\otimes \psi^*) \circ [d, \iota(\beta)]$$avec $\iota(\beta)=\sum\limits_{1\leq i\leq 2n} \beta_i \iota(\hat{e_i})$ et $\iota(\hat{e_i})$ la dérivation de l'algèbre extérieure $\bigwedge (\g^*\oplus \g^*)$  donnée par la contraction. Comme on a $\partial_t m_t=-m_t \circ V(\beta_t)$ il vient alors $$\partial_t M_t =-(m_t\circ V(\beta_t)) \otimes \psi^*= - (m_t\otimes 1) \circ (1\otimes \psi^*) \circ [d, \iota (\beta_t)]=-M_t \circ [d, \iota (\beta_t)] .$$
L'égalité $$\partial_t M_t= - M_t \circ  [d, \iota(\beta_t)] $$ peut se comprendre comme une troisième reformulation de la conjecture de Kashiwara--Vergne dans le complexe des cochaînes.\\

 Cette égalité montre que l'application de Duflo se prolonge en un morphisme d'algèbres de $H(\g, S[\g])$ dans $H(\g, U(\g))$,  précisant ainsi l'isomorphisme démontré dans \cite{Kont}. On retrouve le résultat de \cite{PT}.

\section{Appendice A }

La formule de Kontsevich pour la quantification formelle des variétés de Poisson, montre que l'on peut donner une
expression pour la série  de Campbell-Hausdorff  en utilisant tous les  crochets possibles  à la différence de la formule de Dynkin (cf. \S~\ref{dynkin}).
\subsection{Quantification de Kontsevich}
Dans cette section, afin de faciliter la compréhension de la déformation utilisée dans \cite{AM},  on va rappeler brièvement  la construction de Kontsevich pour la quantification formelle des variétés de Poisson dans le cas du dual des algèbres de Lie.
\subsubsection{Variétés de configurations}
On note $C_{n,m}$ l'espace des configurations de $n$ points distincts
dans le
demi-plan de Poincaré (points de première espèce ou points aériens ) et $m$ points  distincts sur la droite réelle
(ce sont les points de seconde espèce ou points terrestres ), modulo l'action
du groupe $az+b$ (pour $a \in \R^{+*}, b\in \R$).  Dans son article \cite{Kont} Kontsevich construit des
compactifications de ces variétés notées $\overline{C}_{n,m}$. Ce
sont des variétés à coins de dimension $2n-2+m$. Ces variétés ne sont pas connexes pour $m\geq 2$.
On notera par  $\overline{C}^{+}_{n,2}$ la composante  qui contient les configurations où les
 points terrestres sont ordonnés dans l'ordre croissant (ie. on  a $\ov{1}< \ov{2}<\cdots< \ov{m}$).

\begin{figure}[h!]
\begin{center}
\includegraphics[width=6cm]{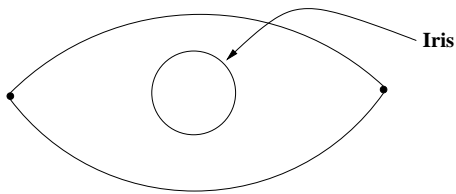}
\caption{\footnotesize La variété $\overline{C}_{2, 0}$.}\label{oeil.eps}
\end{center}
\end{figure}
\subsubsection{Graphes et graphes géométriques}

On note  par $G_{n,2}$ l'ensemble des graphes
étiquetés\footnote{Par graphe étiqueté
on entend un graphe $\Gamma$ muni d'un ordre total sur l'ensemble
$E_\Gamma$ de ses ar\^etes, compatible avec l'ordre des sommets.}  et orientés (les arêtes sont orientées) ayant $n$ sommets de première espèce numérotés $1,2,
\cdots,  n$
et deux sommets de deuxième espèce $\overline{1}, \overline{2}$,  tels que~:

\smallskip
i-  Les arêtes partent des sommets de première espèce. De chaque sommet de première espèce partent exactement deux arêtes.

\smallskip
ii-  Le but d'une ar\^ete est différent de sa source (il n'y a pas
de boucle).

\smallskip
iii- Il n'y a pas d'ar\^ete multiple. \\

\begin{figure}[!h]
\begin{center}
\includegraphics[width=9cm]{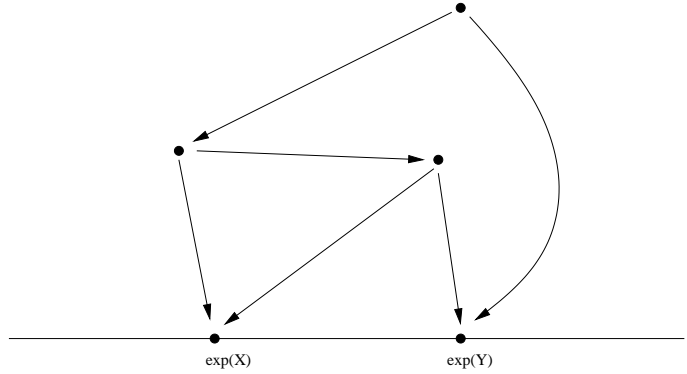}
\caption{\footnotesize Graphe simple de type Lie et de symbole $\G(X, Y)=[[X, [X, Y]], Y]$.}\label{Lie.eps}
\end{center}
\end{figure}

\begin{figure}[h!]
\begin{center}
\includegraphics[width=8cm]{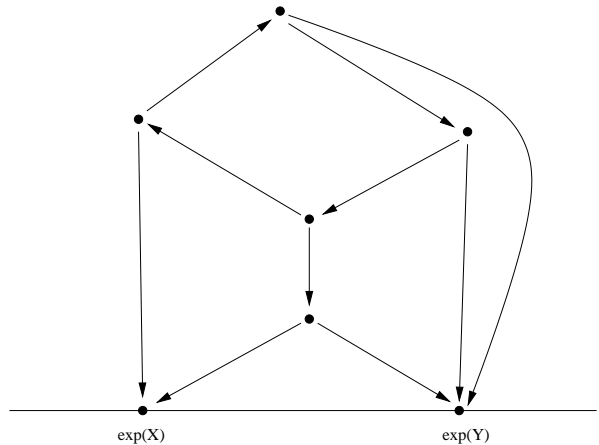}\caption{\footnotesize
  Graphe de type roue  et de symbole $ \Gamma(X,Y)=\tr_\g( \ad X \ad[X,Y]\ad Y\ad
Y).$}\label{roue.eps}
\end{center}
\end{figure}

 Dans le cas linéaire qui nous intéresse, les graphes qui interviennent de manière non triviale (on dira essentiels), sont tels que les sommets de première espèce ne peuvent recevoir qu'au plus une arête. Il en résulte que tout graphe essentiel est superposition de graphes simples de type Lie (graphe ayant une seule racine comme dans Fig.~\ref{Lie.eps} ) ou de type roue (cf. Fig.~\ref{roue.eps} pour un exemple).\\

Les graphes simples essentiels  de type Lie  n'ont pas de symétries.
Par conséquent les graphes de $G_{n,2}$ étiquetés associés à un
graphe géométrique de type Lie (graphe orienté associé pour lequel
on oublie l'étiquetage) sont au nombre de $n!2^n$.

Les graphes
simples de type roue peuvent admettre des
symétries. On notera $m_\Gamma$ le cardinal du groupe de symétries de  $\Gamma$.
 \subsubsection{Fonction d'angle et coefficients}

Soient deux points distincts $(p,q)$ dans le demi-plan de Poincaré muni
de la métrique de Lobachevsky. On note
\begin{equation}
\phi_h(p,q)=Arg\left( \frac{q-p}{q-\overline{p}}\right)
\end{equation}
la fonction d'angle de $C_{2,0}$ dans $\mathbb{S}^1$.
Cette fonction  d'angle s'étend en une fonction régulière à la compactification
 $\overline{C}_{2,0}$.

Si $\Gamma$ est un graphe dans $G_{n,2}$, alors  toute arête $e$ définit par restriction
une fonction d'angle notée $\phi_{e}$ sur la variété $\overline{C}^+_{n,2}$. On note $E_{\Gamma}$
l'ensemble des arêtes du graphe $\Gamma$. Le produit ordonné
\begin{equation}
\Omega_{\Gamma}=\bigwedge _{e \in E_{\Gamma}}\d\phi_{e}
\end{equation}
est donc une $2n$-forme sur  $\overline{C}^+_{n,2}$ variété compacte de dimension $2n$.
\begin{defi} Le poids associé à un graphe $\G$ est par définition

\begin{equation}
w_{\Gamma}=\frac{1}{(2\pi)^{2n}}\int_{\overline{C}^+_{n,2}} \Omega_{\Gamma}.
\end{equation}
\end{defi}

\subsection{Nouvelle formule de Campbell-Hausdorff}\label{sectionKVnew}
Si $\G$ est un graphe  simple de type Lie, on notera  $\Gamma(X, Y)$ le mot dans l'algèbre de Lie libre associée (cf. Fig.~\ref{Lie.eps} pour un exemple). Plus généralement si $\G$ est simple de type roue alors $\G(X, Y)$ sera une fonction de trace (cf. Fig.~\ref{roue.eps} pour un exemple).

\begin{theo}[\cite{Ka}, \cite{AST}]\label{theoKVnew} La série  de Campbell-Hausdorff peut s'écrire en termes de graphes sous la forme d'une série convergente au voisinage de $(0, 0)$
\begin{equation} Z(X,Y)=X+Y +\sum\limits_{n\geq
1}\sum\limits_{\substack{
 \Gamma  \;\mathrm{simple}\\\mathrm{g\acute{e}om\acute{e}trique}\\
\mathrm{de \;type\; Lie}\; (n,2)}}w_{\Gamma} \Gamma(X,Y).
\end{equation}La fonction de densité s'écrit  alors comme série convergente au voisinage de $(0, 0)$ \begin{equation}
D(X,Y)=\exp\Big(\sum\limits_{n\geq 1}\sum\limits_{\substack{
\Gamma \;\mathrm{simple}\\\mathrm{g\acute{e}om\acute{e}trique}\\
\mathrm{de\; type\; Roue \; (n,2)}}} \frac{ w_{\Gamma}
}{m_\Gamma}\G(X,Y)\Big).
\end{equation}
\end{theo}
\subsection{Déformation de Kontsevich}\label{sectiondeformationkont}

On construit  une déformation $2$-dimensionnelle
 de la série  de Campbell-Hausdorff en déformant les coefficients  via un paramètre $\xi \in \overline{C}_{2, 0}$. Pour $\G\in G_{n,2}$ notons $\overline{C}_{\xi}$ la  pré-image de $\xi$ dans l'espace de configurations $\overline{C}_{n+2, 0}$ et
  $$w_\Gamma(\xi)=\frac{1}{(2\pi)^{2n}}\int_{\overline{C}_{\xi}} \Omega_\Gamma.$$
On définit  les déformations $Z_\xi (X,Y) $ et $D_\xi(X, Y)$   en remplaçant dans les formules du théorème \ref{theoKVnew} le coefficient $w_\Gamma$ par sa déformation $w_\Gamma(\xi)$. Ces déformations sont  régulières et admettent  comme conditions aux limites pour $\xi=(0,1)$\footnote{Cette position sur l'axe réel correspond à un coin de $\overline{C}_{2, 0}$.}
$$Z_{(0,1)}(X, Y)=Z(X, Y) \quad\mathrm{et}\quad  D_{(0,1)}(X, Y)=D(X, Y),$$
et pour $\xi=\alpha$ une position sur l'iris\footnote{Cela correspond à une concentration des deux  points  selon un angle $\alpha$.} (cf. Fig~\ref{oeil.eps})$$Z_\alpha(X, Y)=X+Y \quad \mathrm{et}\quad D_\alpha(X, Y)=1.$$
Cette déformation est contrôlée par des équa\-tions
différen\-tielles provenant de l'action tangente du groupe $G$ c'est à dire les champs adjoints qui interviennent dans le conjecture de Kashiwara--Vergne. Ce contrôle provient essentiellement de la formule de Stokes comme  utilisée  dans \cite{Kont}.
\begin{theo}[\cite{To}]\label{theoT} Il existe des séries de Lie  universelles $F_{\xi}(X,Y)$ et $G_{\xi}(X,Y)$
explicites construites en termes de diagrammes, convergentes dans un voisinage de $(0,0)$  et qui sont des $1$-formes régulières sur $\overline{C}_{2,0}$ telles que l'on ait :
\begin{equation}
 \mathrm{d}_\xi Z_\xi(X, Y)= [X,
F_{\xi}(X,Y)]\cdot\partial_{X}Z_\xi(X, Y) + [Y,
G_{\xi}(X,Y)]\cdot\partial_{Y} Z_\xi(X, Y)\end{equation}et
\begin{multline}
\d _\xi D_\xi(X, Y)=\Big([X, F_\xi(X,Y)]\cdot \partial_X +[Y,
G_\xi(X,Y)]\cdot
\partial_Y\Big)D_\xi(X, Y) +\\
\Big(\tr_\g\big(\partial_X F_\xi (X,Y) \circ \ad X+ \partial_Y G_\xi (X,Y)\circ \ad
Y\big) \Big) D_\xi(X, Y).
\end{multline}
\end{theo}

\begin{rema} Pour $\xi$ générique la déformation $Z_\xi$  ne définit pas une loi associative.  Par exemple, la déformation le long de la paupière de $\overline{C}_{2, 0}$ fait intervenir des polynômes de Bernoulli (cf. \cite{To}).
\end{rema}

\begin{rema} On peut montrer \cite{AT} que  la connexion $\gamma_\xi : = (F_\xi, G_\xi) $  est plate,
c'est à dire que l'on a $\d_\xi \gamma_\xi +\frac 12 [\gamma_\xi, \gamma_\xi]_{Lie}=0,$ pour le crochet défini
\S\ref{crochetLie}.
\end{rema}

\begin{figure}[!h]
\begin{center}
\includegraphics[width=8cm]{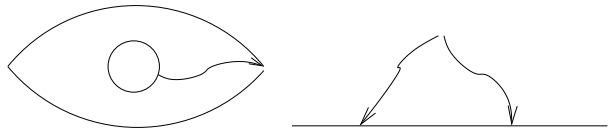}
\caption{\footnotesize Chemin de l'iris jusqu'au coin}\label{chemin.eps}
\end{center}
\end{figure}

\begin{rema} Grâce à la quantification de Kontsevich on construit les déformations  vérifiant les équations (\ref{KV1bis}) et (\ref{KV2bis}). En suivant un chemin comme dans Fig.~\ref{chemin.eps} de l'iris jusqu'au coin, on définit donc des déformations $Z_u(X, Y), D_u(X, Y)$ pour $u\in [0, 1]$ et une fonction $\gamma_u=(F_u, G_u)$ définie dans un voisinage
de $(0, 0)$ à valeurs dans $\g\times\g$.
\end{rema}

\begin{rema} On peut étendre toutes ces constructions au cas de la série de Campbell-Hausdorff avec $n$ arguments $Z(X_1, X_2, \ldots, X_n)$ en considérant un paramètre de déformation dans $\overline{C}_{n, 0}$. On a encore un contrôle  par des équations différentielles comme dans le théorème~\ref{theoT}. La méthode de Alekseev-Meinrenken s'étend sans problème et résout le problème de Kashiwara--Vergne avec $n$ arguments posé dans \cite{bur}.
\end{rema}

\section{Appendice B}
On regroupe dans cet appendice quelques résultats complémentaires sur la formule de Dynkin et on précise certains calculs utiles pour la compréhension de ce texte.
\subsection{La formule de Dynkin}\label{dynkin}
Il existe de nombreuses fa\c cons d'écrire la série  de Campbell-Hausdorff. On
 peut
notamment écrire les développements que l'on obtient en calculant la
dérivée de l'application exponentielle puis en intégrant à nouveau.
Nous allons ici rappeler une autre formule due à   Dynkin. \\

Pour simplifier, on note comme ci-dessous  les crochets successifs
normalisés :
$$[X_1, \ldots, X_n]_*=\frac 1n[X_1,[X_2, \ldots ,[X_{n-1},
X_n]]\ldots]$$et
$$[X_1^{r_1}, \ldots, X_n^{r_n}]_*=[\underset{r_1}{\underbrace{X_1,
\ldots, X_1}}, \ldots, \underset{r_n}{\underbrace{X_n\ldots
X_n}}]_*$$On obtient alors la célèbre formule de Dynkin.

\begin{prop} On a la formule
\begin{multline}
Z(X,Y)=X+Y +\sum_{\substack{m\geq 2}}
\frac{(-1)^{m-1}}{m}\sum_{\substack{p_i+q_i>0}}
\frac{[X^{p_1},Y^{q_1},\ldots ,X^{p_m},Y^{q_m} ]_*}{p_1!q_1!\ldots
p_m!q_m!}.
\end{multline}
\end{prop}
Rappelons comment on obtient cette formule: notons $L_{X, Y}$
l'algèbre de Lie libre engendrée par $X,Y$ \footnote{On appellera
aussi  polynôme
 de Lie en $X, Y$ ou  élément de type Lie, tout élément de $L_{X, Y}$.}
  et $\mathrm{Ass}_{X, Y}$
  l'algèbre associative libre engendrée par $X, Y$.
Pla\c cons nous dans
 l'algèbre enveloppante $U(L_{X, Y})$ qui  rappelons-le s'identifie
  à $\mathrm{Ass}_{X, Y}$ (cf. \cite{Serre}).
 On calcule formellement
$$e^Xe^Y= \sum_{\substack{p, q\geq 0 }}
\frac{X^p}{p!}\frac{Y^q}{q!},$$puis en utilisant le développement de
\[\ln z=\sum_{m\geq 1} \frac{(-1)^{m-1}}{m}(z-1)^m,\] on trouve
formellement
\begin{equation}\nonumber Z(X,Y)= \sum_{\substack{m\geq 1}}
\frac{(-1)^{m-1}}{m}\sum_{\substack{p_i+q_i \geq
1}}\frac{X^{p_1}Y^{q_1}\ldots X^{p_m}Y^{q_m} }{p_1!q_1!\ldots
p_m!q_m!}. \end{equation} On utilise alors une caractérisation des
éléments de l'algèbre de Lie libre dans l'algèbre associative libre due à Dynkin (\cite{Serre} \S
4.4):
un élément $$a=\sum_\alpha  c_\alpha X_{\alpha_1} X_{\alpha_2}\ldots X_{\alpha_n}$$ d'ordre $n$ est dans $L_{X, Y}$ si
et seulement si $$a=\sum_\alpha c_\alpha [X_{\alpha_1}, X_{\alpha_2}, \ldots, X_{\alpha_n}]_*.$$
\begin{rema} La formule de Dynkin n'utilise que des crochets itérés.
\end{rema}

\subsection{La différentielle de l'application exponentielle}\label{sectiondiffexp}

A partir de la formule  de Dynkin on peut mener un calcul explicite
sur les termes à l'ordre $1$ en $y$ (voir \cite{pos}, page 103) ce
qui permet de retrouver les nombres de Bernoulli $b_n$. Ce calcul
n'est pas évident, mais il est faisable.  Rappelons que la série de
Bernoulli est donnée par
\begin{equation}
\sum\limits_{n\geq 0}
\frac{b_nx^n}{n!}=\frac{x}{e^x-1}=1-\frac{x}{2} +\frac{x^2}{12}
-\frac{x^4}{720}+\frac{x^6}{30240}\cdots.
\end{equation}
Les $b_n$ pour $n\geq 3$ impair sont nuls.\\

 En calculant la
série $Z(X,Y)$ à l'ordre $1$ en $Y$, on déduit la formule bien
connue suivante:
\begin{eqnarray}\nonumber
Z(X,Y)\equiv X+Y+ \frac{1}{2} [X,Y] + \frac{1}{12}[X,[X,Y]] +\cdots
\pmod {Y^2}\\\nonumber \equiv X+ \frac{\ad X}{1-e^{-\ad X}}\cdot Y
\pmod {Y^2}.
\end{eqnarray}
On en tire la formule
\begin{eqnarray}\nonumber
\exp_{\g}(X)\exp_{\g}\left(\frac{1-e^{-\ad X}}{\ad X}\cdot Y\right)\equiv
\exp_{\g}(\, X+\, Y+\hspace{-0,3cm}\mod{Y^2}),
\end{eqnarray} et on conclut que  la différentielle de
l'application exponentielle s'identifie à l'endomorphisme
\begin{equation}\label{formulediffexp}Y\mapsto \frac{1-e^{-\ad X}}{\ad X}\cdot Y,\end{equation}
si on utilise la multiplication à gauche\footnote{Si on utilisait la multiplication à droite on trouverait $Y\mapsto \frac{e^{\ad X}-1}{\ad X}\cdot Y$.}
pour identifier $\g$ avec l'espace tangent en $\exp_\g(X)$.

\subsection{Equivalence entre les équations (\ref{KV1}) et (\ref{KV1bis}) pour $Z_t(X,Y)$}\label{equivalence}

Dans cette sous-section nous expliquons comment on passe de
l'équation (\ref{KV1}) à l'équation (\ref{KV1bis}). On suit les références \cite{KV} et \cite{Rou86}.\\

\noi On fait un calcul à l'ordre $1$ en $\epsilon$, comme dans une
dérivée pour \[Z_t(X+ \epsilon [X, F_t], Y+ \epsilon [Y, G_t]).\]

\noi D'après la formule de la différentielle
(proposition~\ref{formulediffexp}) pour l'application exponentielle
on a :

 \begin{multline}\exp_\g(tX+\epsilon [tX,
F_t])=\exp_\g(tX)\exp_\g\Big(\epsilon \frac{1-e^{-\ad t X}}{\ad t X}
[t X, F_t]\Big)=\\\exp_\g(tX)\exp_\g\Big(\epsilon (1-e^{-\ad t X})
F_t\Big).\end{multline}

\noi De même on a

\[\exp_\g(tY+\epsilon [tY, G_t])=\exp_\g\Big(\epsilon (e^{\ad t  Y}-1)
G_t\Big)\exp_\g(tY).\] On en déduit alors

\begin{multline}\label{trace4}\exp_\g\big( t Z_t(X+ \epsilon [X, F_t], Y+ \epsilon [Y, G_t])\big)=\exp_\g(tX+\epsilon [tX, F_t])\exp_\g(tY+\epsilon [tY,
G_t])=\\\exp_\g(tX)\exp_\g\Big(\epsilon (1-e^{-\ad tX}) F_t+
\epsilon (e^{\ad t Y}-1) G_t\Big) \exp_\g(tY).
\end{multline}

\noi L'équation (\ref{KV1bis}) dit que l'on a$$Z_t(X+ \epsilon
[X, F_t], Y+ \epsilon [Y, G_t])=Z_t(X,Y) + \epsilon \partial_t
Z_t(X,Y)=Z_{t+\epsilon}(X,Y).$$Dans ce cas, en utilisant la formule
\[\exp_\g(tY)Z_t(X,Y)\exp_\g(-tY)=Z_t(Y,X),\]le membre de gauche de
(\ref{trace4}) s'écrit
\begin{multline}\label{trace5}
\exp_\g
(tZ_{t+\epsilon}(X,Y))=\exp_\g
((t+\epsilon)Z_{t+\epsilon}(X,Y))\exp_\g(-\epsilon Z_t(X,Y))=\\
\exp_\g((t+\epsilon)X)\exp_\g((t+\epsilon)Y)\exp_\g(-\epsilon
Z_t(X,Y))=\\\exp_\g(tX)\exp_\g\Big(\epsilon X +\epsilon Y -\epsilon
Z_t(Y,X)\Big)\exp_\g(t Y).
\end{multline}

\noi En comparant avec le membre de droite de (\ref{trace4}) on se
trouve que l'on a à l'ordre $1$ en~$\epsilon$ :
\[X+Y-Z_t(Y, X)=(1-e^{-\ad tX}) F_t(X,Y)+ (e^{\ad t Y}-1) G_t(X,Y)\] qui est
exactement l'équation (\ref{KV1}).\\

\noi Réciproquement si l'équation (\ref{KV1}) est vérifiée alors  on aura

\[(1-e^{-\ad tX}) F_t+  (e^{\ad t Y}-1) G_t=X+Y
-\frac 1t\log\left(\exp_\g(tY)\exp_\g(tX)\right).\]En rempla\c cant ce terme
dans (\ref{trace4}) on  retrouve avec le terme de droite de
(\ref{trace5}) ce qui permet de remonter le calcul et de conclure
que l'on a

\[Z_{t+\epsilon}(X,Y)=Z_t\big(X+ \epsilon [X, F_t], Y+ \epsilon [Y,
G_t]\big), \]c'est à dire
\[\partial_t Z_t(X,Y)= \big([X, F_t]\cdot \partial_X+ [Y,
G_t]\cdot \partial_Y\big)Z_t(X,Y),\] qui est bien l'équation
(\ref{KV1bis}).


\begin{thebibliography}{666}

\bibitem[AM00]{AM0} A.~ALEKSEEV, E.~MEINRENKEN --
\textit{The non-commutative Weil algebra} Invent. Math. \textbf{139}  (2000), no.~3, 135--172.

\bibitem[AM02]{AM2} A.~ALEKSEEV, E.~MEINRENKEN --
\textit{Poisson geometry and the Kashiwara--Vergne conjecture.}  C. R. Acad. Sci. Paris,  Sér. I Math.  \textbf{335 } (2002),  no. 9, 723--728.

\bibitem[AM05]{AM5} A.~ALEKSEEV, E.~MEINRENKEN --
\textit{Lie theory and the Chern-Weil homomorphism.} Ann. Sci. École Norm. Sup. (4) \textbf{38}  (2005),  no.~2, 303--338.

\bibitem[AM06]{AM}  A.~ALEKSEEV, E.~MEINRENKEN --  \textit{On the Kashiwara--Vergne conjecture.}
 Invent. Math.  \textbf{164}  (2006), no.~3, 615--634.

 \bibitem[AP]{AP}  A.~ALEKSEEV, E.~PETRACCI --  \textit{On the Kashiwara--Vergne conjecture.}
Journal of Lie Theory \textbf{16}  (2006), 531--538. .

\bibitem[AT]{AT}  A.~ALEKSEEV, C.~TOROSSIAN --  \textit{Star-product for Lie algebras and KV conjecture.}
Preprint 2007.


\bibitem[ADS]{ADS}
M.~ANDLER, A.~DVORSKY, S.~SAHI --
\textit{Kontsevich quantization and invariant distributions on Lie groups}.
 Ann. Sci. École Norm. Sup. (4) \textbf{35} (2002), no.~3, 371--390.

\bibitem[AST]{AST} M.~ANDLER, S.~SAHI, C.~TOROSSIAN -- \textit{Convolution of invariant distributions:
 proof of the Kashiwara--Vergne conjecture. } math.QA/0104100. Lett. Math. Phys. \textbf{ 69}  (2004), 177--203.

\bibitem[B]{bur} E.~BURGUNDER --\textit{Eulerian idempotent and Kashiwara--Vergne conjecture.} math.QA/0612548.




\bibitem[Du70]{Duflo70} M.~DUFLO --
\textit{Caractères des groupes et des algèbres de Lie résolubles.} Ann. Sci. École Norm. Sup. (4) 3 1970 23--74.

\bibitem[Du77]{Duflo77} M.~DUFLO -- \textit{Opérateurs différentiels bi-invariants sur un
groupe de Lie. } Ann. Sci. École Norm. Sup. \textbf{10} (1977),
107-144.



\bibitem[Dix]{Dix} J.~DIXMIER -- \textit{Sur l'algèbre enveloppante d'une algèbre de Lie nilpotente.}
Arch. Math. \textbf{10} 1959 321--326.


\bibitem[Gu]{gutt}
S.~GUTT -- \emph{An explicit {$\sp{\ast} $}-product on the cotangent bundle of a
  {L}ie group}. Lett. Math. Phys. \textbf{7} (1983), no.~3, 249--258.

\bibitem[HC]{HC} Harish-Chandra --
\textit{On some applications of the universal enveloping algebra of a semisimple Lie algebra.}
Trans. Amer. Math. Soc. \textbf{70}, (1951). 28--96.

\bibitem[KV]{KV} M.~KASHIWARA, M.~VERGNE --  \textit{The Campbell-Hausdorff
formula and invariant hyperfunctions. } Invent. Math.
\textbf{47} (1978), 249--272.


\bibitem[Ka]{Ka} V.~KATHOTIA -- \textit{Kontsevich's universal formula for deformation
  quantization and the {C}ampbell-{B}aker-{H}ausdorff formula}. Internat. J.
  Math. \textbf{11} (2000), no.~4, 523--551.

\bibitem[Kly]{kly}
A.~A.~KLJA{\v{C}}KO --  \emph{Lie elements in a tensor algebra}. Sibirsk. Mat. \v
  Z. \textbf{15} (1974), 1296--1304, 1430.


\bibitem[Ko]{Kont}
M.~KONTSEVICH --\textit{Deformation quantization of {P}oisson
manifolds, {I}}.
math.QA/9709040. \textit{Lett. Math. Phys.} \textbf{66
}(2003), no. 3, 157--216.


\bibitem[Mo]{Mo}
T.~MOCHIZUKI --\emph{On the morphism of {D}uflo-{K}irillov type}. J. Geom. Phys.  \textbf{41} (2002), no.~1-2, 73--113.



\bibitem[PT]{PT}
M.~PEVZNER, C.~TOROSSIAN -- \emph{Isomorphisme de Duflo et
cohomologie tangentielle.}  J. Geom. Phys.
\textbf{51} (2004), no.~2,  486--505.

\bibitem[Po]{pos}
M.~POSTNIKOV -- \emph{Le\c cons de géométrie- Groupes et algèbres de Lie}. Ed. Mir, Moscow, 1985.

\bibitem[Rou81]{rou81} F.~ROUVIÈRE -- \textit{Démonstration de la conjecture de
Kashiwara--Vergne pour~$\mathrm{SL}_2(\R)$.} C. R. Acad. Sci. Paris,
\textbf{292} (1981), 657--660.


\bibitem[Rou86]{Rou86} F.~ROUVIÈRE -- \textit{Espaces symétriques et méthode de Kashiwara--Vergne}.  Ann. Sci. école Norm. Sup. (4)  \textbf{19}  (1986),  no. 4, 553--581.


\bibitem[Se]{Serre} J.-P.~SERRE --  \textit{Lie algebras and Lie
  groups}.  New
York-Amsterdam: W. A. Benjamin, Inc. (1965).

\bibitem[Sh]{shoi}
B.~SHOIKHET --  \emph{Tsygan formality and {D}uflo formula}. Math. Res. Lett.
  \textbf{10} (2003), no.~5-6, 763--775.

\bibitem[To]{To} C.~TOROSSIAN -- \textit{Sur la conjecture
combinatoire de Kashiwara--Vergne.}  J. Lie Theory
\textbf{12} (2002),  no. 2, 597--616.

\bibitem[Va]{Va} V.S.~VARADARAJAN -- \textit{Lie Groups, Lie Algebras and their Representations}. Graduate Texts in Mathematics,  Ed Springer-Verlag, 1984.

\bibitem[Ve]{Ve}  M.~VERGNE --  \textit{ Le centre de l'algèbre enveloppante et la formule de
 Campbell-Hausdorff.} C. R. Acad. Sci. Paris,  Sér. I Math.  \textbf{329 } (1999),  no. 9, 767--772.

\end{thebibliography}
\end{document}